\documentclass[final,12p,times,authoryear]{elsarticle}
\usepackage[margin=1.5cm]{geometry}

\usepackage[utf8]{inputenc}
\usepackage{amssymb}
\usepackage{amsmath}
\usepackage{bm}
\usepackage{lipsum}
\usepackage{subcaption}
\usepackage{float}
\setcitestyle{square}

\usepackage{pgfplots}                   
\pgfplotsset{width=7cm, compat=newest}
\usetikzlibrary{patterns}
\tikzset{every label/.style={font=\footnotesize,inner sep=1pt}}
\newcommand{\stencilpt}[4][]{\node[circle,fill = red,draw,inner sep=1.5pt,label={above:#4},#1] at (#2) (#3) {}}

\usepackage{algpseudocode}
\usepackage{algorithm}

\usepackage{tikz}
\usepackage{hyperref}
\usepackage{ulem}
\usepackage{wrapfig}
\usepackage{standalone}
\usetikzlibrary{cd}
\usetikzlibrary{positioning}
\usetikzlibrary{fit}
\usetikzlibrary{calc,patterns,angles,quotes}
\tikzset{set/.style={draw,circle,inner sep=0pt,align=center}}

\definecolor{morange}{RGB}{255,127,14}
\definecolor{mblue}{RGB}{31,119,180}
\definecolor{mred}{RGB}{214,39,40}
\definecolor{mpurple}{RGB}{148,103,189}
\definecolor{mgreen}{RGB}{44,160,44}

\definecolor{ccccccc}{RGB}{204,204,204}
\definecolor{cffffff}{RGB}{255,255,255}

\journal{Computers \& Fluids}

\begin{document}

\begin{frontmatter}

\title{Data-Driven Adaptive Gradient Recovery for Unstructured Finite Volume Computations}

\author[first,second]{G. de Rom\'{e}mont}
\author[first]{F. Renac}
\author[second]{F. Chinesta}
\author[first]{J. Nunez}
\author[first]{D. Gueyffier}
\affiliation[first]{organization={ONERA/DAAA},
            addressline={29 Av. de la Division Leclerc}, 
            city={Châtillon},
            postcode={92320}, 
            state={France}}

\affiliation[second]{organization={ENSAM},
            addressline={151 Bd de l'Hôpital}, 
            city={Paris},
            postcode={75013}, 
            state={France}}


\begin{abstract}
We present a novel data-driven approach for enhancing gradient reconstruction in unstructured finite volume methods for hyperbolic conservation laws, specifically for the 2D Euler equations. Our approach extends previous structured-grid methodologies to unstructured meshes through a modified DeepONet architecture that incorporates local geometry in the neural network. The architecture employs local mesh topology to ensure rotation invariance, while also ensuring first-order constraint on the learned operator.
The training methodology incorporates physics-informed regularization through entropy penalization, total variation diminishing penalization, and parameter regularization to ensure physically consistent solutions, particularly in shock-dominated regions. The model is trained on high-fidelity datasets solutions derived from sine waves and randomized piecewise constant initial conditions with periodic boundary conditions, enabling robust generalization to complex flow configurations or geometries.
\\
Validation test cases from the literature, including challenging geometry configuration, demonstrates substantial improvements in accuracy compared to traditional second-order finite volume schemes. The method achieves gains of 20-60\% in solution accuracy while enhancing computational efficiency. A convergence study has been conveyed and reveal improved mesh convergence rates compared to the conventional solver.
\\
The proposed algorithm is faster and more accurate than the traditional second-order finite volume solver, enabling high-fidelity simulations on coarser grids while preserving the stability and conservation properties essential for hyperbolic conservation laws. This work is a part of a new  generation of solvers that are built by combining Machine-Learning (ML) tools with traditional numerical schemes, all while ensuring physical constraint on the results. 

\end{abstract}


\end{frontmatter}



\let\thefootnote\relax\footnotetext{* corresponding author : guillaume(dot)romemont(at)protonmail(dot)com (G. de Rom\'{e}mont)}

\section{Introduction}\label{introduction}
Computational Fluid Dynamics (CFD) commonly relies on nonlinear hyperbolic partial differential equations (PDEs) to model a wide range of complex fluid behaviors. Even when the initial or boundary conditions are smooth, these equations can lead to the formation of discontinuities in finite time \cite[Sec. 2.4.2]{toro2000centred}. 
In this scope, standard unstructured second order finite volume methods have become essential in CFD allowing the discretization of complex geometries like aircraft wings or turbine blades. However, achieving accurate solutions requires extremely fine grid resolution, particularly in regions with steep gradients, boundary layers, or turbulent flow, where large spatial or temporal scales dramatically increase computational cost and memory requirements. The need for fine grids stems from the fact that gradient reconstruction accuracy directly impacts solution quality. Machine learning presents a transformative opportunity to break this cycle by learning optimal reconstruction strategies that can maintain high accuracy even on relatively coarse meshes.
\par 
In this context, a number of works assessed the use of machine learning in computational physics, offering new avenues to augment or even replace traditional physics-driven numerical models. These data-driven strategies have proven effective in terms of accelerating numerical simulations, enhancing precision or robustness. In scenarios where the data maintains temporal but mainly spatial regularity, methods originating from image and video analysis such as multilayer perceptrons (MLPs), convolutional neural networks (CNNs), U-nets, recurrent neural networks (RNNs), long short term memory (LSTM) can be leveraged for scientific computing, including computational physics and fluid dynamics, e.g. forecasting physical fields, identifying model parameters, generating high-resolution outputs...
\\
However, the reliance of these techniques on structured data such as Cartesian grids, constant time intervals, or uniform input sizes restricts their effectiveness when applied to unstructured or irregular datasets. This presents a notable obstacle in domains where complex geometry is commonplace within simulations, particularly in high-resolution, physics-based simulations.
\par 
To overcome these limitations, growing efforts are being directed toward the development of methods capable of handling unstructured data, common in many real-world problems. Many methods are derived from Reduce Order Modeling (ROM), like Proper Orthogonal Decomposition (POD) or Singular Value Decomposition (SVD) and have been applied to many areas successfully with unstructured meshes like ocean models \cite{fang2009pod}, turbulent flows \cite{berkooz1993proper} , compressible aerodynamics \cite{bui2003proper, manzoni2015reduced}. Others methods are developed mapping unstructured or sparse data onto a structured grid allowing the use of structured-based machine learning techniques like CNNs \cite{liu2015sparse, coscia2023continuous}. An other promised and logical direction is the use of Graph Neural Networks (GNNs) where they excel by their capability to handle adaptive geometries. GNNs can be divided into three categories: convolutional \cite{hamilton2017inductive}, attentional \cite{velickovic2017graph} and message passing \cite{gilmer2017neural}. GNNs have been used to model various complex simulations \cite{sanchez2020learning}.
\par 
However, hyperbolic systems develop discontinuous solutions and standard neural networks surrogates struggle with shocks. For example, continuous activation functions are needed for proper learning during backpropagation  (the process of reconstructing the gradient of an optimization function through Automatic Differentiation). A significant advancement has been the ability to precisely satisfy certain physical or numerical constraints by incorporating learned models into a fixed equation of motion. Indeed, hyperbolic conservation laws must preserve certain integral quantities (mass, momentum, energy). Unconstrained ML models may violate these conservation properties, leading to nonphysical solutions. 
Constraints can be enforced through conservative network architectures, Lagrange multipliers in the loss function, or post-processing correction steps. 
\par 
A range of methods have been developed for meshless solutions with Lagrange multipliers for physical guidance with strong formulation like PINNs \citep{raissi2018hidden, raissi2019physics} or weak formulation \citep{yu2018deep, kharazmi2019variational, bar2019learning} more adequate for hyperbolic PDEs. But Lagrange multipliers only define a physical guidance and not a hard constraint on the solver, meaning the physics are only approximated. 
\par 

To this end, several hybrid physics-informed machine learning approaches have been developed to integrate neural networks within established numerical frameworks for solving forward problems . More specifically in finite volume solvers particularly suitable for hyperbolic PDEs \citep{jessica2023finite, li2023finite, ranade2021discretizationnet, stevens2020finitenet}. In this paradigm, different works assessed the use of ML for improving shocks capturing methods \citep{stevens2020enhancement,magiera2020constraint}, enhancing flux limiters \citep{nguyen2022machine, schwarz2023reinforcement}, enhancing corrections with artificial viscosity \citep{bruno2022fc} or the flux itself \citep{bezgin2021data,morand2024deep}. Some works focused their methodologies for unstructured data \citep{cen2024deep, morand2024deep}. The closest work related to our methodology is the optimization of a compact high-order variational reconstruction on triangular grids \citep{zhou2024machine}. In this study, an artificial neural network is employed to predict the optimal values of the derivative weights on cell interfaces. These interfaces represent the free parameters of the variational reconstruction.
\par
Pursuing a similar goal, \cite{bar2019learning} introduced a learning-based approach to gradient interpolation that matches the accuracy of conventional finite difference schemes while operating on significantly coarser grids. This methodology has been integrated into classical finite volume solvers, extending its application to problems such as passive scalar advection \cite{zhuang2021learned} and the Navier-Stokes equations \cite{kochkov2021machine}. However, these approaches are tailored to specific governing equations, for periodic boundary conditions and are highly unstable. 
\par 
This paper presents a machine learning gradient optimization framework for second-order finite volume schemes on triangular unstructured grids extending and generalizing previous work of \cite{de2024data}. The work specifically targets the 2D Euler equations. Using a supervised learning framework with high resolution solutions for the database, a geometry-aware neural network is trained to predict free parameters in order to correct the gradient used in the scheme. 
The corrections will be explored for two types of gradients, namely the Green-Gauss (GG) and the Least Square (LSQ) gradients. The neural network architecture is inspired from the DeepONet architecture \citep{lu2019deeponet} with local inputs values and geometry. The geometry inputs are angles between each neighbors cells and allows more flexibility concerning the skewness of some elements. The model is trained on a dataset with high-resolution solutions made of random Riemann or the fluxed initial conditions with periodic boundary conditions. Lagrange multipliers such as the Total Variation regularizer and the entropy regularizer have been added to the loss to physically constrain the output solution with key properties inherent to hyperbolic conservations laws. The entropy regularizer aims at selecting the physically relevant weak solution \citep[Thm 6.1]{godlewski2013numerical}.
As the method developed is local, the work can be generalized to all types of geometries with different types of boundary conditions. 
Numerical results for two-dimensional Euler test cases demonstrate that the machine learning
optimized gradient finite volume scheme can achieve a gain of 20\% to 60\% in solution accuracy for a same mesh configuration. A computational performance review has also been conveyed against the traditional finite volume solver. 
\\
The challenges addressed by this paper are the same as in \cite{de2024data}, meaning stability, accuracy and computational performance. 
\par
The main innovation of this paper is the derivation of a highly constrained accurate and robust methodology to accurately compute solution for hyperbolic PDEs on unstructured mesh. The neural network can be seen as a subgrid correction operator, thus constrained and allowing super-resolution and physically-consistent solutions.
\par
The remainder of this paper is organized as follows: Section \ref{sec:model_problem} provides an overview of hyperbolic conservation laws with a deeper look at 2D Euler equations. In Section \ref{sec:finite_volume_solver}, we describe the finite volume solver and the implementation of boundary conditions. Section \ref{sec:data_driven} introduces the complete algorithm alongside the Machine Learning model used and several regularization designed to ensure that the neural network produces physically consistent solutions to forward problems, particularly in the presence of strong shocks. In Section \ref{sec:results}, the method is validated against a range of benchmark equations and test cases from the literature, demonstrating strong performance in both accuracy and stability. Finally, Section \ref{sec:Computational} presents numerical experiments evaluating the algorithm’s computational efficiency against the traditional finite volume solver.
%
\section{Model problem}\label{sec:model_problem}

\subsection{Nonlinear hyperbolic conservation laws} 

We are here interested in the approximation of first-order nonlinear hyperbolic conservation laws and consider initial and boundary value problems in $d\geq1$ space dimensions of the form

\begin{subequations}\label{eq:hyperbolic}
\begin{align}
    \partial_t \textbf{w} + \nabla\cdot\textbf{f}(\textbf{w}) &= 0 \quad \text{in }\Omega\times(0,\infty), \label{eq:hyperbolic_a}\\
		\textbf{w}(\bm{x},0) &= \textbf{w}_0(\bm{x}) \quad \text{in }\Omega,  \label{eq:hyperbolic_b}\\
    \bm{B}(\textbf{w},\textbf{w}_{bc},\bm{n}) &= 0 \quad \text{on }\partial\Omega\times(0,\infty),  \label{eq:hyperbolic_c}
\end{align}
\label{eq:hyperbolic_eq}
\end{subequations}

\noindent with $\Omega \subset \mathbb{R}^d$ a bounded domain, $\textbf{w}:\Omega \times [0,T] \rightarrow \Omega^a\subset\mathbb{R}^r$ denotes the vector of $r$ conservative variables with initial data $\textbf{w}_0\in L^\infty(\Omega)$ and $\bm{f}:\Omega^a  \rightarrow \mathbb{R}^r \times \mathbb{R}^d$ are the physical fluxes. The solution is known to lie within a convex set of admissible states $\Omega^a$ \citep[Sec. IV]{dafermos2005hyperbolic}. Boundary conditions are imposed on $\partial\Omega$ through the boundary operator in (\ref{eq:hyperbolic_c}) and some prescribed boundary data $\textbf{w}_{bc}$ defined on $\partial\Omega$, while $\bm{n}$ denotes the unit outward normal vector to $\partial\Omega$. The operator $\bm{B}$ depends on the type of condition to be imposed and the equations under consideration. Examples are provided in section \ref{sec:BC_examples}.

Solutions to (\ref{eq:hyperbolic}) may develop discontinuities in finite time even if $\textbf{w}_0$ is smooth, therefore the equations have to be understood in the sense of distributions. Weak solutions are not necessarily unique and (\ref{eq:hyperbolic}) must be supplemented with further admissibility conditions to select the physical solution. We here focus on entropy inequalities of the form 
\cite[sec. I.5]{godlewski2013numerical}\cite[sec. III]{dafermos2005hyperbolic}\cite[sec. III.8]{leveque1992numerical}

\begin{equation}\label{eq:entropy_ineq}
    \partial_t \eta(\textbf{w}) + \nabla\cdot\bm{q}(\textbf{w})\leq 0,
\end{equation}

\noindent where $\eta:\Omega^a\rightarrow\mathbb{R}$ is a convex entropy function \cite[sec. III.8]{leveque1992numerical} and $\bm{q}(\textbf{w}):\Omega^a\rightarrow\mathbb{R}^d$ is the entropy flux satisfying the compatibility condition 

\begin{equation}
    \nabla_{\textbf{w}} \eta(\textbf{w})^T \times \nabla_{\textbf{w}}\textbf{f}_i(\textbf{w}) = \nabla_{\textbf{w}}\bm{q}_i(\textbf{w})^T, \quad 1\leq i\leq d.
\end{equation}

The inequality (\ref{eq:entropy_ineq}) becomes an equality for smooth solutions, while strict convexity of $\eta(\textbf{w})$ implies hyperbolicity of (\ref{eq:hyperbolic_a}) \citep{godlewski2013numerical,dafermos2005hyperbolic}.

We here aim at satisfying these properties at the discrete level with our data-driven finite volume scheme. Satisfying these entropy conditions ensure uniqueness for a weak solution by determining which shocks are physically feasible \citep{dafermos2005hyperbolic}.

\subsubsection{Compressible Euler equations}

The work of this article will focus on the compressible Euler equations in two space dimensions, $d=2$, for which

\begin{align*}
    \textbf{w} &= (\rho, \rho \bm{v}^T, E)^T, \\
    \bm{f}(\textbf{w}) &= \big(\rho \bm{v},\rho \bm{v}\bm{v}^T + p{\bf I}_d, (E+p)\bm{v}\big)^T,
\end{align*}

\noindent where ${\bf I}_d$ is the identity matrix of size $d$, while $\rho$, $\bm{v}$ and $E$ denote the density, velocity vector and total energy, respectively. We close the system with an equation of state for a polytropic ideal gas law, so  $E=\tfrac{p}{\gamma - 1} + \tfrac{1}{2}\rho |\bm{v}|^2$ with $\gamma>1$ the ratio of specific heats and $p$ the pressure. Note that it will be convenient to consider the numerical discretization as a function of the primitive variables 

\begin{equation*}
    \textbf{u} = (\rho, \bm{v}^T, p)^T.
\end{equation*}
Note that for scalar equations, $u = w$.
\par
This system is hyperbolic in every unit direction $\bm{n}$ over the set of admissible states $\Omega^a=\{\textbf{w}\in\mathbb{R}^{d+2}:\,\rho>0,\bm{v}\in\mathbb{R}^d,E-\tfrac{1}{2}\rho |\bm{v}|^2>0\}$ with eigenvalues $\bm{v}\cdot\bm{n}-c$, $\bm{v}\cdot\bm{n}$, $\bm{v}\cdot\bm{n}+ c$, where $c=\big(\tfrac{\gamma p}{\rho}\big)^{1/2}$ is the speed of sound. Even if not unique, a typical entropy-flux pair $(\eta, \mathbf{q})$ for the Euler equations is given by

\begin{equation}\label{eq:entropy}
    \eta = -\rho s, \quad
    \mathbf{q}(\textbf{w}) = -\rho s \bm{v},
\end{equation}

\noindent where $s = C_v \log\left(\tfrac{P}{\rho^\gamma}\right)$ is the physical specific entropy defined by the Gibbs relation
\begin{equation*}
    Tds = \dfrac{1}{\gamma - 1}d\left( \dfrac{p}{\rho}\right) - \dfrac{p}{\rho^2}d\rho
\end{equation*}

\subsection{Boundary conditions}\label{sec:BC_examples}

To a given boundary $\Gamma\subset\partial\Omega$ in (\ref{eq:hyperbolic_c}), we associate an admissible boundary state

\begin{equation}\label{eq:boundary_state}
 \textbf{w}_\Gamma:\Omega^a\times\mathbb{S}^{d-1}\ni(\textbf{w},\bm{n})\mapsto\textbf{w}_\Gamma(\textbf{w},\bm{n})\in\Omega^a,
\end{equation}

\noindent which we assume to be admissible, $\textbf{w}_\Gamma(\Omega^a\times\mathbb{S}^{d-1})\subset\Omega^a$, and consistent: ${\bf B}(\textbf{w},\textbf{w}_{bc},\bm{n})=0$ implies $\textbf{w}_\Gamma(\textbf{w},\bm{n})=\textbf{w}$.

Considering a unit normal vector $\bm{n} = [n_x, n_y]^T$ pointing outwards the domain. 
The flux Jacobian $\bm{A}(\textbf{w},\bm{n}) = \partial_{\textbf{w}} \bm{f}(\textbf{w})\cdot\bm{n}$ has eigenvalues $(\lambda_i)_{1\leq i\leq r}$ and is diagonalized as $A(\textbf{w},  \bm{n}) = R^{-1} \Lambda R$ with $\Lambda = \text{diag}(\lambda_i)$. Depending on the sign of the eigenvalues, we introduce the upwind decomposition of $A(\textbf{w},  \bm{n})$ as follows
\begin{equation*}
    A(\textbf{w},  \bm{n}) = A^+(\textbf{w},  \bm{n}) + A^-(\textbf{w},  \bm{n}), 
\end{equation*}

\noindent with 
\begin{equation*}    
    \quad A^\pm(\textbf{w},  \bm{n}) = R^{-1} \Lambda^\pm  R, \quad \Lambda^\pm=\dfrac{1}{2}\text{diag}\big(\lambda_i\pm|\lambda_i|\big).
\end{equation*}
\par 
Each eigenvalue is associated to a given characteristics. The characteristics are used to determine the number of boundary conditions to be imposed \citep{goncalvesdasilva:cel-00556980,hartmann2002adaptive}, which correspond to the inflow characteristics that propagate from outside to inside the computational domain.
\par 

For the compressible Euler equations, the flux Jacobian $\bm{A}(\textbf{w},\bm{n}) = \partial_{\textbf{w}} \bm{f}(\textbf{w})\cdot\bm{n}$ has eigenvalues

\begin{equation*}
    \lambda_1=\bm{v}\cdot\bm{n}-c,\quad \lambda_2=\dots=\lambda_{d+1}=\bm{v}.\bm{n}, \quad \lambda_{d+2}=\bm{v}\cdot\bm{n} + c.
\end{equation*}
Namely, $c = \sqrt{\dfrac{\gamma p }{\rho}}$ is the speed of sound.

\paragraph{Farfield conditions} We can apply characteristic boundary conditions on the farfield boundary from a freestream state $\textbf{w}_{\infty}$ with 
\begin{equation}
    B(\textbf{w}, \textbf{w}_{bc}, \bm{n}) = A^- (\textbf{w} - \textbf{w}_{\infty}),
\end{equation}
\noindent where $\textbf{w}_{\infty}$ denotes a given freestream state. For instance, for a supersonic inflow boundary condition, we have $A^-=A$, corresponding to the Dirichlet boundary conditions
\begin{equation}
    \textbf{w}_{\Gamma}(\textbf{w}, \bm{n}) =  \textbf{w}_{\infty},
\end{equation}
while for the supersonic outflow boundary condition, we have $A^+=A$ and $A^-=0$, corresponding to an extrapolation condition
\begin{equation}
    \textbf{w}_{\Gamma}(\textbf{w}, \bm{n}) =  \textbf{w}.
\end{equation}

\paragraph{Slip boundary condition} Considering an impermeability condition, $\bm{v}.\bm{n} = 0$, at a wall $\Gamma_w \subset \partial \Omega$, the associated boundary data is $\textbf{w}_\Gamma(\textbf{w}, \bm{n}) = (\rho, \rho (\bm{v}- (\bm{v}.\bm{n})\bm{n})^T, \rho E)^T$. The condition is commonly imposed  through the use of a mirror state $2 \textbf{w}_\Gamma(\textbf{w}, \bm{n}) - \textbf{w} = (\rho, \rho (\bm{v}- 2(\bm{v}.\bm{n})\bm{n})^T, \rho E)^T$. The mirror state $\bm{u}^+_\Gamma(\bm{u}^-, \bm{n})$ follows from imposing a linear reconstruction interpolating the right and left states at the interface : $\frac{1}{2}(\textbf{w}^- + \textbf{w}^+_\Gamma(\textbf{w}^-, \bm{n})) = \textbf{w}_\Gamma(\textbf{w}^-, \bm{n})$ hence $\textbf{w}^+_\Gamma(\textbf{w}^-, \bm{n})) = 2 \textbf{w}^-_\Gamma(\textbf{w}, \bm{n}) - \textbf{w}^-$

\paragraph{Periodic condition} 
This boundary is established when physical geometry of interest and expected flow pattern are of a periodically repeating nature. This reduces computational effort in our problems.

%
%
\section{Finite volume solver} \label{sec:finite_volume_solver}

We describe below the finite volume solver used to generate the data. We use a MUSCL reconstruction of the slopes \citep{van1979towards} on 2D unstructured mesh to get a formally second-order scheme. In section \ref{finite_volume_solver}, we introduce the finite volume method in 2D, while the treatment of boundary conditions is described in section \ref{sec:FV_BCs}.
\par
The mesh will be generated using triangles as a partition of the space domain $\Omega$ into a set of $N_{\Omega}$ disjoint cells $C_i$:
\begin{equation*}
    \Omega = \displaystyle \bigcup_{i=1}^{N_{\Omega}}C_i
\end{equation*}
All cells areas $|C_i|$ are non-zeros and the mesh is untangled. We here consider a cell-centered finite volume scheme.
\par
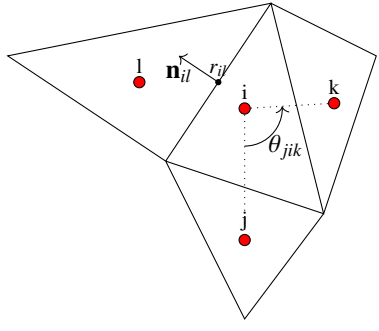
\begin{figure}
    \centering
        \begin{tikzpicture}[scale=0.7]
        \coordinate (i) at (1.5,1);
        \coordinate (j) at (1.5,-1.5);
        \coordinate (k) at (3.2,1.1);
        \coordinate (l) at (-0.5,1.5);
        \draw (0,0) -- (3,-1) -- (2,3) -- (0,0);
        \draw (3,-1) -- (1.5,-3) -- (0,0);
        \draw (3,-1) -- (4,2) -- (2,3);
        \draw (2,3) -- (-3, 2) -- (0,0);
        \stencilpt{1.5,1}{i}{i};
        \stencilpt{1.5,-1.5}{j}{j};
        \stencilpt{3.2,1.1}{k}{k};
        \stencilpt{-0.5,1.5}{l}{l};
         \draw[->] (1,1.5) -- (0.25,2) node[anchor=north]{$\mathbf{n}_{il}$};
         \draw[dotted] (i) -- (j) ;
         \draw[dotted] (i) -- (k) ;
         \node[draw,circle,fill = black, scale=0.1, label=$r_{il}$] at (1,1.5) {$r_{il}$};
        \pic [draw, ->, "$\theta_{jik}$", angle eccentricity=1.5] {angle = j--i--k};
    \end{tikzpicture}\\
    \caption{Unstructured mesh in 2D The red bullets denote the center of the cells. $\mathbf{n}_{il}$ is the unit normal to the face $S_{il}$ with center $r_{ij}$. $\theta_{jik}$ is the measured angle between the points $j$, $i$ and $k$.}
    \label{fig:unstructured_mesh} 
\end{figure}

Some geometric notations pictured on Figure \ref{fig:unstructured_mesh}:
\begin{equation*}
\begin{cases}    
S_{ij} &= \partial \overline{C}_i \cap \partial \overline{C}_j = \text{ common face between } C_i \text{ and } C_j\\
S_{ib} &= \partial \overline{C}_i \cap \partial \overline{\Omega} = \text{ face of } C_i \text{ on boundary of } \Omega\\
\mathcal{N}_i &= \{C_j:\partial \overline{C}_i \cap \partial \overline{C}_j \ne \emptyset\}\\
\textbf{n}_{ij} &= \text{outward normal of face } S_{ij} \text{ from } C_i \text{ to } C_j  \\
r_i& = (x_i, y_i), \text{ the cell center }\\
r_{ij} &= \text{ center of face } S_{ij}
\end{cases}
\end{equation*}
On Figure \ref{fig:unstructured_mesh}, $\mathcal{N}_i = \{C_l, C_k, C_j\}$. 
\noindent
For the sake of simplicity, instead of writing $C_j \in \mathcal{N}_i$, we will simply say that $j \in \mathcal{N}_i$
\subsection{Finite volume solver} \label{finite_volume_solver}
Integrating equation (\ref{eq:hyperbolic_a}) over a cell $C_i$ gives:
\begin{equation}
    \displaystyle \int_{C_i}  \partial_t \textbf{w}(\textbf{x},t) dV + \displaystyle \int_{C_i} \nabla\cdot\textbf{f}(\textbf{w}) dV = 0
\end{equation}
We define the cell average value as 
\begin{equation*}
    \overline{\textbf{w}}_{i}(t) = \dfrac{1}{|C_i|}\int_{C_i} \textbf{w}(\textbf{x},t) dV
\end{equation*}
And using the divergence theorem, we obtain
\begin{equation}
    |C_i|\partial_t \overline{\textbf{w}}_{i} + \displaystyle  \int_{\partial C_i} \textbf{f}(\textbf{w})\cdot \textbf{n} dS = 0
\end{equation}
\begin{equation}
    |C_i|\partial_t \overline{\textbf{w}}_{i} + \displaystyle \sum_{j\in\mathcal{N}_i}\int_{S_{ij}} \textbf{f}(\textbf{w})\cdot \textbf{n} dS + \displaystyle \sum_{S_{ib}\in\partial\Omega}\int_{S_{ib}} \textbf{f}(\textbf{w})\cdot \textbf{n} dS = 0
\end{equation}
 We have to approximate the flux integral by quadrature. For first order and second order accurate schemes, it is enough to use mid-point rule of integration.
 \begin{equation}
     \int_{S_{ij}} \textbf{f}(\textbf{w})\cdot \textbf{n} dS \simeq \left[\textbf{f}(\textbf{w})\cdot \textbf{n}\right]_{ij}|S_{ij}|
 \end{equation}
 Giving two states $\textbf{w}_{ij}$ and $\textbf{w}_{ji}$ coming from cells $C_i$ and $C_j$, we will use a numerical flux function of Godunov-type
  \begin{align*}
      \left[\textbf{f}(\textbf{w})\cdot \textbf{n}\right]_{ij} &\simeq H(\textbf{w}_{ij}, \textbf{w}_{ji},\textbf{n}_{ij}), \\
      \left[\textbf{f}(\textbf{w})\cdot \textbf{n}\right]_{ib} &\simeq H(\textbf{w}_{ib}, \textbf{w}_{b},\textbf{n}_{ib}),
 \end{align*}
 thus giving
 \begin{equation}
    |C_i|\partial_t \overline{\textbf{w}}_{i} + \displaystyle \sum_{j\in\mathcal{N}_i} H(\textbf{w}_{ij}, \textbf{w}_{ji},\textbf{n}_{ij}) + \displaystyle \sum_{S_{ib}\in\partial\Omega} H(\textbf{w}_{ib}, \textbf{w}_{b},\textbf{n}_{ib}) = 0
\end{equation}

As a reference solver, we use a second-order MUSCL \citep{van1979towards} finite volume scheme with a Rusanov numerical flux \citep{lax2005weak, rusanov1961calculation} and a Venkatakrishnan limiter \citep{venkatakrishnan1995convergence} which is used because it is a smooth differentiable function as opposed for example to the minmod slope limiter \citep{roe1986characteristic}. The flux can be written
\begin{align*}
    H(\textbf{w}_{ij}, \textbf{w}_{ji},\textbf{n}_{ij}) &= \frac{1}{2}\left[\textbf{f}(\textbf{w}_{ij}) + \textbf{f}(\textbf{w}_{ji})\right]\cdot\mathbf{n}_{ij} - \frac{1}{2}\rho(\textbf{w}_{ij}, \textbf{w}_{ji})(\textbf{w}_{ij} - \textbf{w}_{ji})
\end{align*}
\par
To achieve  a second order accuracy, the primitives variables are reconstructed inside each cell using the limiter.
\begin{align*}
    \textbf{u}_{ij} &=\textbf{u}_i + \phi_i (r_{ij} - r_{i})\cdot\nabla \textbf{u}_i \\
    \textbf{u}_{ji} &=\textbf{u}_j + \phi_j(r_{ij} - r_{j})\cdot\nabla \textbf{u}_j
\end{align*}
then $\textbf{w}_{ij}$ and $\textbf{w}_{ji}$ are evaluated from $\textbf{u}_{ij}$ and $\textbf{u}_{ji}$.
\par 
\par 
Time integration is conveyed with the traditional explicit Euler method. For faster training, this method is preferred to the more accurate high order Runge-Kutta methods. Note that even if training is made with explicit Euler, inference can be made with other temporal integrations schemes. The time-step $\Delta t$ is chosen constant regarding 
\begin{equation}
    Co = \dfrac{\Delta t}{\underset{j}{\operatorname{min}}(\sqrt{|C_j|})} =constant
    \label{eq:Co_definition}
\end{equation} 
where the classical CFL condition is defined as
\begin{equation}
    C_{FL} = Co \times |\lambda|_{max}
    \label{eq:CFL_definition}
\end{equation} 
with $|\lambda|_{max}$ the maximum eigenvalue of the flux Jacobian $\bm{A}$.
\par
Indeed, as seen in Sec. \ref{sec:data_driven}, we will train our neural network with high resolution solutions and choosing this type of CFL-like parameter allows us to make sure that two solutions are at the same given time $t$ for accurate comparison. If there is an overshoot for the velocity for the ML solution, using the standard CFL condition will be become useless as the time step will differ.
\\
During training, $Co = 0.03$ for stability reasons and during inference $Co$ is usually chosen of the same order but can be increased or decreased. 

\subsection{Gradient computation} \label{sec:gradient}
The gradient is computed using either the Green-Gauss method denoted $\nabla_{GG}$
\begin{equation*}
   \frac{1}{|C_i|} \displaystyle\int_{C_i}\nabla \textbf{u} dV = \frac{1}{|C_i|}\displaystyle\int_{\partial C_i}\textbf{u} \textbf{n} dS 
\end{equation*}
approximated by
\begin{equation}
   \nabla_{GG}\textbf{u}_i = \frac{1}{|C_i|}\displaystyle \sum_{j\in \mathbf{N}_i}\frac{1}{2}(\textbf{u}_i + \textbf{u}_j)\textbf{n}_{ij}|S_{ij}|
   \label{eq:GG_def}
\end{equation}
or the Least Square (LSQ) method with the notation $\nabla_{LSQ}$.
To compute $\nabla \textbf{u}_i = (\textbf{a},\textbf{b})$, we define $\Delta \textbf{u}_j = \textbf{u}_j - \textbf{u}_i$, $\Delta x_j = x_j - x_i$ and $\Delta y_j = y_j - y_i$ giving 
\begin{equation}
   \nabla_{LSQ}\textbf{u}_i = 
   \begin{pmatrix}
   \sum_j \omega_j \Delta x_j^2 & \sum_j \omega_j \Delta x_j \Delta y_j\\
   \sum_j \omega_j \Delta x_j \Delta y_j & \sum_j \omega_j \Delta y_j^2
   \end{pmatrix}^{-1}
   \begin{pmatrix}
        \sum_j \omega_j \Delta x_j \Delta \textbf{u}_j \\
         \sum_j \omega_j \Delta y_j \Delta \textbf{u}_j 
   \end{pmatrix}
   \label{eq:LSQ_def}
\end{equation}
with the weight function $\omega_i$ taken of the form
\begin{equation*}
    \omega_j = \dfrac{1}{|r_j-r_i|^2}
\end{equation*}

The least squares method is a reliable approach for obtaining accurate gradient estimates. However, in highly anisotropic grids, this approach can result in unstable schemes. The implementation of distance-based weighting allows to improve robustness and accuracy of the reconstruction by lowering the influence of distant cells.

\subsection{Venkatakrishnan limiter} \label{sec:limiter}
This limiter is a smooth modification of the min-max limiter \cite{venkatakrishnan1995convergence}. Smoothness of the limiter is here useful for improving the convergence during the learning process based on gradient-based descent algorithm.
We first define
\begin{equation*}
    \begin{cases}
        \textbf{u}_i^m &= \underset{j\in\mathcal{N}_i}{\operatorname{min}}(\textbf{u}_j, \textbf{u}_i) \\
        \textbf{u}_i^M &= \underset{j\in\mathcal{N}_i}{\operatorname{max}}(\textbf{u}_j, \textbf{u}_i)
    \end{cases}
\end{equation*}

and 
\begin{equation*}
    \Delta_{ij} = (r_{ij} - r_i)\cdot\nabla \textbf{u}_j
\end{equation*}
such that the limiter is defined as
\begin{equation}
  \phi_{ij} =
    \begin{cases}
      L(\textbf{u}_i^M - \textbf{u}_i,\Delta_{ij}) & \text{if $\Delta_{ij}>0$ }\\
      L(\textbf{u}_i^m - \textbf{u}_i,\Delta_{ij}) & \text{if $\Delta_{ij}<0$ }\\
      0 & \text{otherwise}
    \end{cases}       
\end{equation}
with $L(a,b) = \dfrac{a^2+2ab+\omega}{a^2+2b^2+ab}$.

\subsection{Boundary conditions}\label{sec:FV_BCs}
 The implementation of far-field boundary conditions for a flow problem is contingent upon two prerequisites. Primarily, the truncation of the domain must not have a discernible impact on the flow solution when compared to that of an infinite domain. Secondly, any outgoing disturbances must not be reflected back into the flow field. The boundary conditions are usually defined from prescribed freestream values.
 \par
 In order to achieve this, the computational domain is modified by the introduction of so-called "ghost cells" situated outside of the boundary. The introduction of fictitious flow in the ghost cells will yield the desired boundary conditions at the edge.

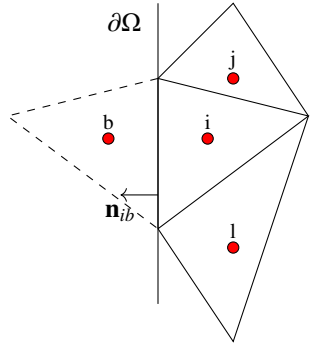
\begin{figure}[H]
    \centering
        \begin{tikzpicture}
        \draw[dashed]  (0,2) -- (-2,1.5) -- (0,0);
        \draw (0,0) -- (0,2) -- (2,1.5) -- (0,0);
        \draw (2,1.5) -- (1,3) -- (0,2);
        \draw (0,0) -- (1,-1.5) -- (2,1.5);
        \stencilpt{-0.66,1.2}{b}{b};
        \stencilpt{0.66,1.2}{i}{i};
        \stencilpt{1,2.}{j}{j};
        \stencilpt{1,-0.25}{l}{l};
        \filldraw[black] (-0.1,3) node[anchor=north east]{$\partial\Omega$};
        \draw (0,-1) -- (0,3);
         \draw[->] (0,0.45) -- (-0.5,0.45) node[anchor=north]{$\mathbf{n}_{ib}$};
    \end{tikzpicture}
    \caption{Ghost cell approach at boundaries, the ghost cell is represented in dashed line}
    \label{fig:boundary}
\end{figure}

\paragraph{Supersonic inflow}
The conservative variables on the boundary are determined exclusively on the basis of freestream values. 
\begin{equation*}
    \textbf{w}_{\Gamma} = \textbf{w}_{b}
\end{equation*}
\paragraph{Supersonic outflow} The conservative variables on the boundary are determined exclusively on the basis of inner flow field values. 
\begin{equation*}
    \textbf{w}_{\Gamma} = \textbf{w}_{i}
\end{equation*}
\paragraph{Subsonic inflow} 
Four characteristic variables are prescribed based on the freestream values, which are defined as follows: One characteristic variable is derived from the interior of the physical domain. This results in the following set of boundary conditions:
\begin{align*}
    p_{\Gamma} &= \frac{1}{2}\left[ p_{b}+p_{i}-\rho_{0}c_{0}\bigl[(\bm{v}_{b}- \bm{v}_{i})\cdot \bm{n}\bigr] \right] \\
    \rho_{\Gamma} &= \rho_{b}+\frac{(p_{\Gamma}-p_{b})}{c_{0}^{2}} \\
    \bm{v}_{\Gamma} &= \bm{v}_{b}-\frac{\bm{n}(p_{b}-p_{\Gamma})}{p_{0}c_{0}} 
\end{align*}
where $p_0$ and $c_0$ represent a reference state. The reference state is normally set equal to the state at the interior.
\paragraph{Subsonic outflow}
\begin{align*}
        p_{\Gamma} &= p_{b}\\
        \rho_{\Gamma}&=\rho_{i}+\frac{\left(p_{\Gamma}-p_{i}\right)}{c_{0}^{2}}\\
        \bm{v}_{\Gamma}&=\bm{v}_{i}-\frac{\bm{n}\left(p_{i}-p_{\Gamma}\right)}{\rho_{0}c_{0}}
\end{align*}
with $p_{b}$ being the prescribed static pressure. 
\paragraph{Slip-wall boundary} The implementation of the slip-wall boundary condition can be applied by removing the normal velocity component.
\begin{align*}
        p_{-1} &= p_{i} \\
        \rho_{-1} &= \rho_{i} \\
        \bm{v}_{-1} &= \bm{v}_{i} -  2(\bm{v}_{i}\cdot \bm{n})\bm{n}
\end{align*}
\paragraph{Periodic condition} 
Due to the periodicity condition, the first ghost-cell layer $\textbf{w}_{-1}$ corresponds to the  inner flow field value at the opposite periodic boundary.

%
%
\section{Data driven solution for hyperbolic equations} \label{sec:data_driven}

\subsection{Learning the derivatives}
We set two kind of derivatives to learn denoted with $\hat{.}$ inspired from the Green-Gauss Eq. (\ref{eq:GG_def}) and least square methods Eq. (\ref{eq:LSQ_def}) respectively
\begin{equation}
    \begin{aligned}
            \hat{\nabla}_{GG}\textbf{u}_i &= \dfrac{1}{|C_i|}\displaystyle \sum_{j\in \mathcal{N}_i}\left(\left(\frac{1}{2}+\alpha_j\right)\textbf{u}_i + \left(\frac{1}{2}-\alpha_j\right)\textbf{u}_j\right)\textbf{n}_{ij}|S_{ij}| \\
             &= \nabla_{GG}\textbf{u}_i + \dfrac{1}{|C_i|}\displaystyle \sum_{j\in \mathcal{N}_i}\alpha_j\left(\textbf{u}_i -\textbf{u}_j\right)\textbf{n}_{ij}|S_{ij}|
    \end{aligned}
    \label{eq:ML_GG}
\end{equation}
 and 
 \begin{equation}
      \begin{aligned}
        \hat{\nabla}_{LSQ}\textbf{u}_i &= 
        \begin{pmatrix}
           \sum_j \omega_j \Delta x_j^2 & \sum_j \omega_j \Delta x_j \Delta y_j\\
           \sum_j \omega_j \Delta x_j \Delta y_j & \sum_j \omega_j \Delta y_j^2
           \end{pmatrix}^{-1}
           \begin{pmatrix}
                \sum_j \omega_j (1+\alpha_j)\Delta x_j \Delta \textbf{u}_j \\
                 \sum_j \omega_j(1+\alpha_j) \Delta y_j \Delta \textbf{u}_j 
           \end{pmatrix}\\
         &= \nabla_{LSQ}\textbf{u}_i + 
         \begin{pmatrix}
           \sum_j \omega_j \Delta x_j^2 & \sum_j \omega_j \Delta x_j \Delta y_j\\
           \sum_j \omega_j \Delta x_j \Delta y_j & \sum_j \omega_j \Delta y_j^2
           \end{pmatrix}^{-1}
           \begin{pmatrix}
                \sum_j \omega_j \alpha_j\Delta x_j \Delta \textbf{u}_j \\
                 \sum_j \omega_j\alpha_j \Delta y_j \Delta \textbf{u}_j 
           \end{pmatrix}
\end{aligned}
\label{eq:ML_LSQ}
\end{equation}

All $\alpha_j$ presented in Eq. (\ref{eq:ML_GG}) and Eq. (\ref{eq:ML_LSQ}) will be local outputs of a neural networks. These outputs can be seen as weighting corrections of the gradient for every node neighbor value. 
 \par
We remark that the computation of the gradient is exact for linear functions using least squares and the derivatives using this approach are in general first order accurate. On smooth and symmetric stencils, we can obtain close to second order accuracy. For this reason, and for the rest of the paper if not specified, the ML-gradient will be computed using the LSQ inspired method.
\par
The training of a neural network inside a classic numerical solver is made possible by writing the entire program in a differentiable programming framework. Indeed the computation of the gradients of the loss is made possible using a differentiable CFD solver. This differentiable framework allows Automatic Differentiation \citep{baydin2018automatic} for any parameter of the solver. Multiple frameworks have been recently developed like TensorFlow \citep{abadi2016tensorflow}, Pytorch \citep{paszke2017automatic}, Flux.jl \citep{innes2018flux}, JAX \citep{bradbury2018jax} and are more and more efficient on hardware accelerator (GPUs, TPUs). These user-friendly frameworks make it easier to incorporate neural networks techniques into scientific codes. We have implemented our solver using TensorFlow Eager \citep{agrawal2019tensorflow}. The entire code allows batching computations with same mesh size examples, training time is thus accelerated using gradient descent batching.

\subsection{A DeepONet architecture to efficiently compute the gradients}
For now, the algorithm for enhancing a standard finite volume solver is presented on Algorithm \ref{alg:training}. The architecture of the model $\mathcal{NN}$ will be defined subsequently.

\begin{algorithm}[H]
\caption{A ML-enhanced finite volume solver algorithm}
\label{alg:training}
    \begin{algorithmic}
        \State \textbf{Inputs :} A neural-network Model $\mathcal{NN}$
        \State Optimized training parameters $\Theta^*$
        \State A second-order finite volume solver $FV(\mathbf{u}, \nabla\mathbf{u)}$
        \State A function $\hat{\nabla}()$ for gradient reconstruction (Eq. (\ref{eq:ML_GG}, \ref{eq:ML_LSQ})
        \State A time-step $\Delta t$, and a number of time steps $N_T$
        \State Initial condition $\mathbf{u}^0$
        \State Discretization of domain $\Omega$
         \State Load model parameters $\Theta^*$ onto model $\mathcal{NN}$
         \For{$t = 0,\Delta t,...,N_T\Delta t$}
              \State $\bm{\alpha}^t = \mathcal{NN}(\mathbf{u^t)}$
              \State $\hat{\nabla}\mathbf{u}^t = \hat{\nabla}(\mathbf{u}^t, \bm{\alpha}^t)$
              \State $\mathbf{u}^{t+\Delta t} = FV(\mathbf{u^t, \hat{\nabla}\mathbf{u}^t)}$
         \EndFor
    \end{algorithmic}
\end{algorithm}
\par
The architecture is derived  from the DeepONet \cite{lu2019deeponet}. The DeepONet architecture is inspired by the universal approximation theorem \cite{hornik1989multilayer, cybenko1989approximation} to learn any nonlinear operator $\mathcal{G}$ accurately and efficiently from a relatively small dataset. $\mathcal{G}$ is an operator taking an input function $\mathbf{u}$. The DeepONet as presented in \citep{lu2019deeponet} is composed of two neural networks: the branch net that encodes the input function at a fixed number of sensors $\mathbf{x}$ and the trunk net that encodes the sensors. Mathematically, it is given as
\begin{equation}
    \mathcal{G}(\mathbf{u})(\mathbf{x}) \approx \sum_{k=1}^{p} b_{k}(\mathbf{u}) t_{k}(\mathbf{x})
\end{equation}
The trunk and branch networks outputs respectively $p$ scalars $t_{k}$ and $b_{k}$.
\begin{figure*}
    \centering
    \includegraphics[width=\linewidth]{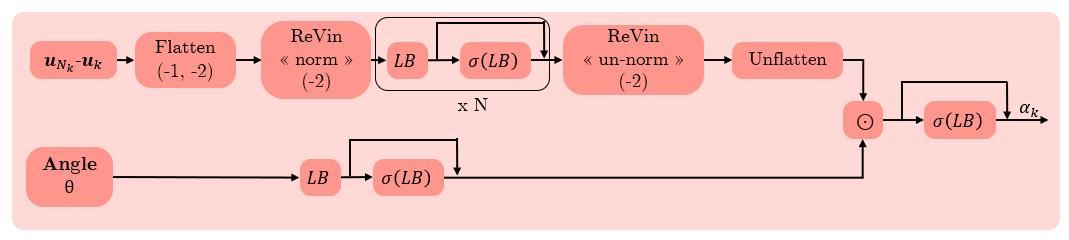}
    \caption{ML architecture inspired from the DeepONet architecture \cite{lu2019deeponet}. Namely, $LB$ is a linear block, $\sigma$ an activation function and $\odot$ is the scalar product between the branch net and the trunk net. The model outputs a local vector $\alpha_k$ necessary for Eq. (\ref{eq:ML_GG}) or (\ref{eq:ML_LSQ}).}
    \label{fig:architecture}
\end{figure*}
\par
Combining a modified version of the DeepONet architecture presented on Figure \ref{fig:architecture} and Eq. (\ref{eq:ML_GG}) or (\ref{eq:ML_LSQ}), we learn the operator $\mathcal{G}_\Theta$ with  $\Theta$ the learnable parameters
\begin{equation}
    \mathcal{G}_\Theta(\mathbf{u})(\mathbf{x}) \approx \nabla \mathbf{u}
\end{equation}
\par
First of all, as a derivative is a local operator, the learned operator needs to be also local and will take only local inputs or geometries.
The input function is $\mathbf{u}_{\mathcal{N}_i} - \mathbf{u}_i$ and the sensors are angles at $\theta_{jik}$, $\theta_{kil}$ and $\theta_{lij}$. Here, $\theta_{jik} = \measuredangle r_j r_i r_k$ represents the measured angle between the vectors $\overrightarrow{ r_i r_j} $ and $\overrightarrow{ r_i r_k}$ represented for example on Figure \ref{fig:unstructured_mesh}. $\mathbf{u}_{\mathcal{N}} - \mathbf{u}$ is the solution vector for each primitive variables of the difference between the neighbors values and the value at cell center. Taking for example the geometry represented in  Figure \ref{fig:unstructured_mesh} for a variable $y$.
\begin{equation*}
    y_{\mathcal{N}_i} - y_i = 
    \begin{pmatrix}
        y_j - y_i\\ y_k - y_i \\y_l - y_i
    \end{pmatrix}
\end{equation*}
This choice of input allows to enforce directly the condition $\hat{\nabla}\textbf{(u+c)} =\hat{\nabla}\textbf{u}$ where $\textbf{c}$ represents a constant. This choice is also enforced on the derivation of Eq. (\ref{eq:ML_GG}) and Eq. (\ref{eq:ML_LSQ}) as the inputs are of the same form in the sums.
\par 
The utilization of angles is based on the fact that it allows the network to understand the distribution of the data locally. There is no use of the distances between cells centers as input for the ML algorithm as the distances or the geometry in general is hard-constrained by the formulation of Eq. (\ref{eq:ML_GG}) or (\ref{eq:ML_LSQ}). We have to note that because of the lack of referential, i.e. all angles are computed in reference to one-another, the ML architecture proposed is rotation-equivariant working well in our case. Indeed, all direction needed for the computation of the gradient is computed in the formulation of Eq. (\ref{eq:ML_GG}) or (\ref{eq:ML_LSQ}).
\par
The flatten and un-flatten layers allows us to establish a link between the different variables $(\rho, \mathbf{v}, p)$, it improves generalization and robustness. Moreover, the use of these type of operation allows us to increase the number of neurons as the matrix used will be bigger because the input vector is bigger, it prevent us to increase the number of trainable parameters by adding ResNets sequentially and increase inference time.
\par
A last point to note is the use of a normalization and un-normalization layer \cite{kim2021reversible} for an improved generalization and robustness for a larger range of values for the variables.

\subsection{Loss} \label{sec:loss}
The training phase adjusts the machine learning parameters to optimize the weights and biases, minimizing the discrepancy between a computationally expensive high-resolution simulation and the model-generated simulation on a coarse grid. This refinement is achieved through supervised learning, where the total loss function $\mathcal{L}_{tot}$ is defined as the cumulative pointwise error between the predicted and reference primitive variables. Since the numerical scheme is inherently stable, employing a flux-limited second-order approach, there is no need to accumulate errors over multiple steps for additional stabilization. The Loss is defined as in \cite{de2024data} and adjusted for unstructured mesh.
\begin{equation}
    \mathcal{L}_{tot} = 100\cdot\dfrac{||\textbf{u}^{ML} - \textbf{u}^{ref}||_2}{|| \textbf{u}^{ref}||_2} + \lambda_{TVD}\mathcal{L}_{TVD} + \lambda_{ent}\mathcal{L}_{ent} + \lambda_{reg}\mathcal{L}_{reg}
\end{equation}
where $\textbf{u}^{ML}$ represents the ML solution, while $\textbf{u}^{ref}$ denotes the reference solution obtained from a fine grid using the reference finite volume scheme and projected back onto the coarse grid.
\par 
Although the method remains stable and robust even without penalization, incorporating penalization provides additional regularization and enhances robustness. This is particularly beneficial around shock regions, where spurious oscillations may still arise despite the presence of a flux limiter and for long term predictions. These additional penalization terms are referred to as soft constraints \citep{raissi2018hidden}. The Entropy and Total Variation Diminishing (TVD) penalizations were first introduced in \cite{patel2022thermodynamically} for hyperbolic PDEs.

\paragraph{Entropy inequality penalization} 
As defined in \cite{de2024data} and following Eq. (\ref{eq:entropy}), the entropy inequality penalization is 
\begin{equation}
    \mathcal{L}_{ent} = \dfrac{1}{N_{\Omega}}\sum_{i\in \Omega} max\left(0, |C_i|(\eta^{t+\Delta t} - \eta^{t}) + \dfrac{|C_i|\Delta t}{2} (\nabla\cdot \textbf{q}^ {t+\Delta t} + \nabla\cdot \textbf{q}^ {t})\right)^2
\end{equation}

\paragraph{Total variation penalization} 
The total variation of a differentiable function $f$ in multiple dimensions is defined as
\begin{equation*}
    TV(f, \Omega) = \int_\Omega |\nabla f|dV
\end{equation*}
so we define the TVD penalization as 
\begin{equation*}
    \mathcal{L}_{TVD} =\displaystyle \sum_{i \in \Omega} max(0,  |\nabla_{LSQ} \textbf{u}_i^{t+\Delta t}| - |\nabla_{LSQ} \textbf{u}_i^{t}|)
\end{equation*}
\par 
This entropy inequality penalization push the model to learn weak unique solution of the hyperbolic PDEs regarding entropy inequality constraint \citep[Thm 3.4]{lax1973hyperbolic}.

\paragraph{Regularization term} 
The regularization term imposes a penalty on the weights and biases of the neural network to mitigate overfitting by preventing excessively large values. This helps to suppress oscillations, particularly in the presence of strong shocks. The penalization is formulated as:
\begin{equation}
    \mathcal{L}_{reg} = \lVert W \rVert_1
\end{equation}
with $W$ being all the parameters of the neural network.

\subsection{Boundary conditions treatment} \label{sec:BC_ML}
As specified in \citep{de2024data}, the local operator learned can only "see" immediate neighborhood, but boundary conditions often encodes relationships between the solution and its environment and can be seen as global constraints. Ghost cells described in Sec. \ref{sec:FV_BCs} can be seen as artificial data creation essentially making up non-physical quantities. 
\par
In the scope of reconstructing accurate and physical gradients, and for good generalization properties, the choice has been made to not train specific ML algorithms for each type of boundary conditions but to set the contribution of the ML algorithm to zero at the boundaries. So that 
\begin{equation*}
    \forall i \in \{ i| \partial \overline{C}_i \cap \partial \overline{\Omega}\neq \emptyset\}, \bm{\alpha}_i = 0
\end{equation*}
We can note that this choice of implementation is particularly useful for slip-wall boundary conditions, this type of boundary alongside the use of ghost cells creates an artificial discontinuity for the velocity.
\\
Furthermore, this choice leads to the classical second order accuracy near boundaries. And because the Euler equations are only considered in this study, there is no need for mesh refinement near boundaries.

\subsection{Dataset} \label{sec:dataset_description}
The dataset is constructed on a fine grid, and subsequently, the solution obtained is projected back onto a coarse grid for the purposes of training and achieving super-resolution properties. The projection operator onto the coarse grid is a area-weighted mean to keep the conservation laws properties.
All solutions utilized in the dataset are computed on a uniform fine mesh, which is derived from a coarse mesh through a standard mesh refinement process. The refinement involves the division of each triangle in the coarse mesh into four triangles. This approach facilitates data batching, thereby enhancing the efficiency of the training process.
\par 
As in \cite{de2024data}, the training data comprises sine waves and randomized piecewise constant initial conditions on a $\Omega = [0,1]^2$ domain that have been derived using pre-defined, parametrized functions. A more exhaustive derivation of the functions used in the dataset is presented in Annex \ref{sec:dataset}. 
\par
For good generalization properties, the initial conditions primitives values are extended to high values ($\in[-6,6]$ for the scaled Euler equations) for both regular or discontinuous initial conditions. It ensures that the model learns to handle the full range of values it will encounter in practice. Furthermore, when models are trained on inputs with high magnitudes, they learn to be less sensitive to the absolute scale of the inputs. This helps the model generalize better to test data that might have different scaling than the training data.

\subsection{Model training} \label{sec:model_training}

The model has been trained on a dataset that exclusively uses periodic boundary conditions. It ensures better generalization performances as periodic boundaries naturally generates a greater number of shocks interaction. These additional shocks contribute to the enrichment of the training data, thereby facilitating the model's learning of more robust patterns and dynamics.
\par
For the dataset, the coarse mesh has about $2614$ cells, and $10456$ cells for the fine mesh used to compute the reference solutions. With this configuration, we generate 8 different solutions for 2000 time-steps with random initial conditions as described in Annex \ref{sec:dataset}, making 16000 solutions frames in the dataset. A validation dataset is also generated with 4 different initial conditions
\par 
In practice and for reduced inference time, we use 2 blocks ($N=2$ on Figure \ref{fig:architecture}) for the primitives and only one block for the angles, making the total number of trainable parameters to 1332. As the ML method is called for each cells and at each time-step, a fast inference time is needed to be computationally efficient.
\par
For the loss function, $\lambda_{TVD}$, $\lambda_{ent}$ and $\lambda_{reg}$ are usually chosen as 1\% to 10\% of the size of the first term of the loss $100\cdot\dfrac{||\textbf{u}^{ML} - \textbf{u}^{ref}||_2}{|| \textbf{u}^{ref}||_2}$,  and they are dependent on choice of the dataset. But in our case, we set $\lambda_{TVD} = 10^{-6}$, $\lambda_{ent}=10^5$ and $\lambda_{reg} = 10^{-4}$. 
\par 
For training, the Lion optimizer \cite{chen2023symbolic} is used with a learning rate of $lr= 6e-5$ and an exponential scheduler \cite{li2019exponential} with the scheduler parameter set to $0.9$. For generally 30 epochs, the training takes approximately 20mn on one RTX-6000 GPU. In the continuation of \cite{de2024data}, there is no need for rolling steps during training time to improve the robustness of the scheme.

%
%

\section{Results} \label{sec:results}
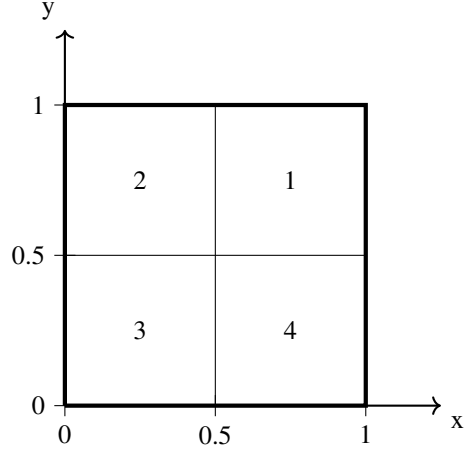
\begin{figure}[H]
    \centering
    \begin{tikzpicture}[scale=4]
        \draw[ultra thick] (0,0) -- (1,0) -- (1,1) -- (0,1) -- (0,0);
        \draw[thick,->] (0,0) -- (1.25,0) node[anchor=north west] {x};
        \draw[thick,->] (0,0) -- (0,1.25) node[anchor=south east] {y};
        \draw (0 cm,1pt) -- (0 cm,-1pt) node[anchor=north] {$0$};
        \draw (0.5 cm,1pt) -- (0.5 cm,-1pt) node[anchor=north] {$0.5$};
        \draw (1 cm,1pt) -- (1 cm,-1pt) node[anchor=north] {$1$};
        \draw (1pt,0 cm) -- (-1pt,0 cm) node[anchor=east] {$0$};
        \draw (1pt,0.5 cm) -- (-1pt,0.5 cm) node[anchor=east] {$0.5$};
        \draw (1pt,1 cm) -- (-1pt,1 cm) node[anchor=east] {$1$};
        
        \draw (0.5,0) -- (0.5,1) ;
        \draw (0,0.5) -- (1,0.5) ;

        \node at (0.25, 0.75) (A) {$2$}; 
        \node at (0.75, 0.75) (A) {$1$}; 
        \node at (0.25, 0.25) (A) {$3$}; 
        \node at (0.75, 0.25) (A) {$4$}; 
    \end{tikzpicture}
    \caption{Division of the domain $\Omega$ for the Riemman test cases}
    \label{fig:Riemman_geometry}
\end{figure}

All test cases presented on a $\Omega = [0,1]^2$ domain are reference cases from \cite{lax1998solution}. The results will focus on Riemman problems dividing the domain $\Omega$ into four quadrants where the initial data are constant in each quadrant as seen on Figure \ref{fig:Riemman_geometry}.

All Riemann problems are posed such that the solutions to the four one-dimensional Riemann problems between quadrants exhibit precisely one propagating disturbance (i.e., a shock, rarefaction, or contact waves). All results presented are computed using the LSQ based gradient method if not specified. The gain is computed as
\begin{equation*}
    \operatorname{gain} = 100\cdot\dfrac{\mathcal{L}_{coarse} - \mathcal{L}_{ML}}{\mathcal{L}_{coarse}}
\end{equation*}
where $\mathcal{L}_{coarse}$ and $\mathcal{L}_{ML}$ is the 1-norm between the fine solution and respectively the coarse solution done by the traditional solver or the ML solution. Mathematically, at a given time $t$ they are computed as 
\begin{align*}
    \mathcal{L}_{coarse} &= \displaystyle \sum^{N_\Omega}_j \lVert \mathbf{u}_j^{ref} - \mathbf{u}_j^{coarse}\rVert_1 \\
    \mathcal{L}_{ML} &= \displaystyle \sum^{N_\Omega}_j \lVert \mathbf{u}_j^{ref} - \mathbf{u}_j^{ML}\rVert_1
\end{align*}
As described in Sec. \ref{sec:dataset_description}, $\mathbf{u}_j^{ref}$ is a high resolution solution projected back onto the coarse mesh using weighted area subsampling. $\mathbf{u}^{coarse}$ and $\mathbf{u}^{ML}$ are the solutions obtained from respectively the traditional solver and the ML-enhanced solver on a coarse mesh. These errors are also named $\mathcal{L}_1$ on the labels of the figures.
\par
For the time-stepping scheme, $Co=0.01$. This can be considered a low value, but it allows a good time-discretization accuracy for a single-step explicit Euler method. Additionally, this ensures that the errors presented are mainly spatial discretizations errors and are not related to the time-discretization. 

\subsection{Riemann test cases}
All Riemman test cases have been computed with $max(|C_i|)\leq3e-4$ representing 5150 cells. Numbering of the test cases is taken from the original paper \citep{lax1998solution}.
\paragraph{Case 6}
The initial data $(\rho, u, v, p)$ for each corresponding quadrants on Figure \ref{fig:Riemman_geometry} are 
\begin{table}[H]
    \centering
    \begin{tabular}{|c|c|}
        \hline
        $(2, 0.75, 0.5, 1)$ & $(1, 0.75, -0.5, 1)$  \\ \hline
        $(2, -0.75, 0.5, 1)$ & $(3, -0.75, -0.5, 1)$  \\ \hline
    \end{tabular}
\end{table}

In this case, four shock waves interact. All boundary conditions have been set to subsonic outflow. One thing to note is the preservation of symmetry and the notably high gain obtained, around $34\%$.

\begin{figure}
     \centering
     \begin{subfigure}[b]{0.35\textwidth}
         \centering
         \includegraphics[width=\textwidth]{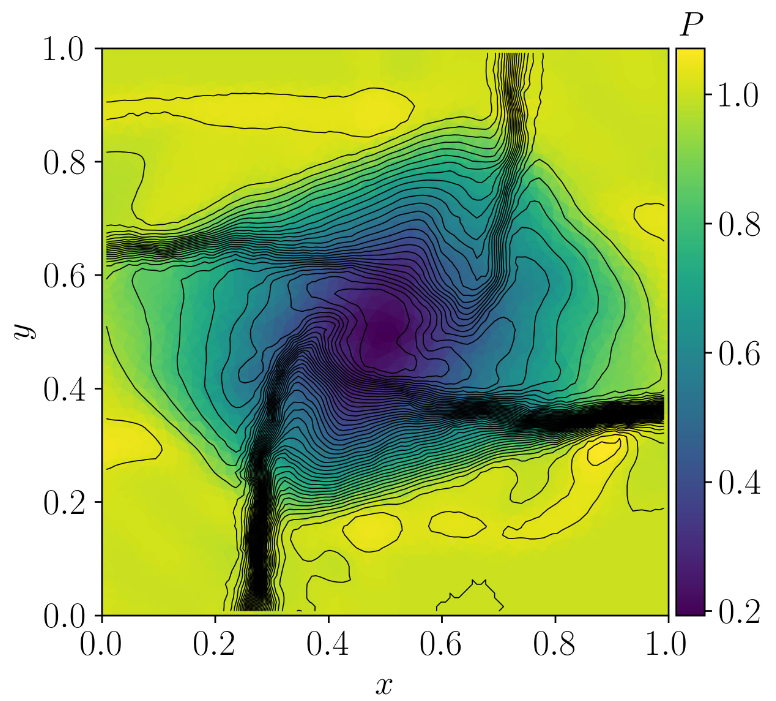}
         \caption{}
         \label{fig:p_contours_test6}
     \end{subfigure}
     \begin{subfigure}[b]{0.4\textwidth}
         \centering
         \includegraphics[width=\textwidth]{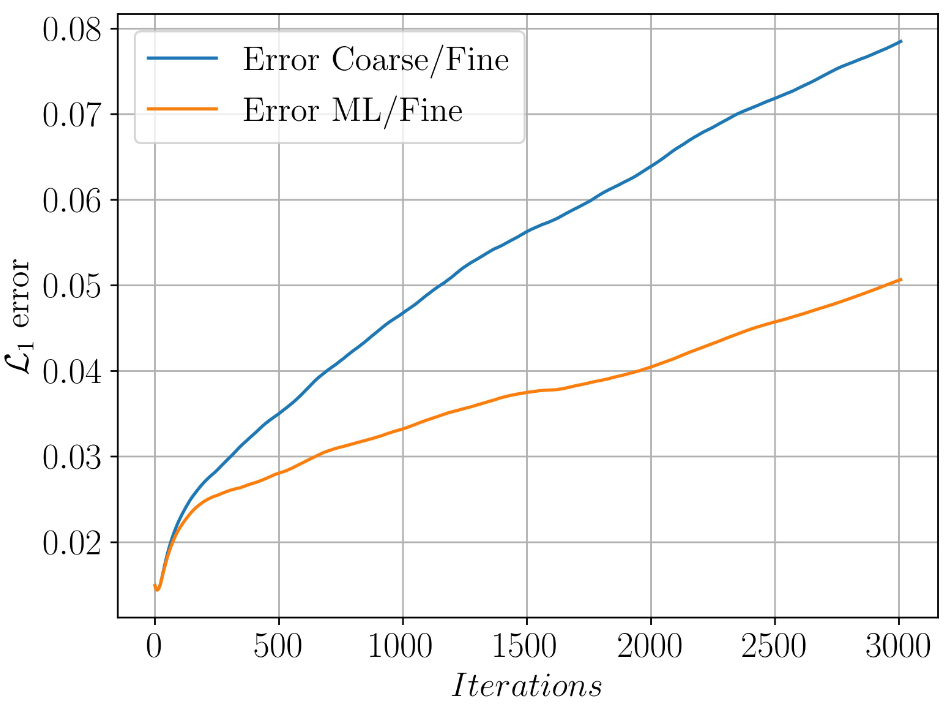}
         \caption{}
         \label{fig:error_test6}
     \end{subfigure}
     \caption{Results for the 2D Riemann problem case 6 \citep{lax1998solution}. Density color map is overlayed by 30 density contours. The computations were performed on the  $(x, y) \in (0, 1) \times (0, 1)$ square with 5150 cells. (a) Pressure $p$ , (b) Error across iterations.}
     \label{fig:test_case_6}
\end{figure}

\paragraph{Case 11}
The initial data $(\rho, u, v, p)$ for each corresponding quadrants on Figure \ref{fig:Riemman_geometry} are 
\begin{table}[H]
    \centering
    \begin{tabular}{|c|c|}
        \hline
        $(0.5313, 0.8276, 0, 0.4)$ & $(1, 0.1, 0.1, 1)$  \\ \hline
        $(0.8, 0.1, 0, 1.4)$ & $(0.5313, 0.1, 0.7276, 0.4)$  \\ \hline
    \end{tabular}
\end{table}

\begin{figure}
     \centering
     \begin{subfigure}[b]{0.35\textwidth}
         \centering
         \includegraphics[width=\textwidth]{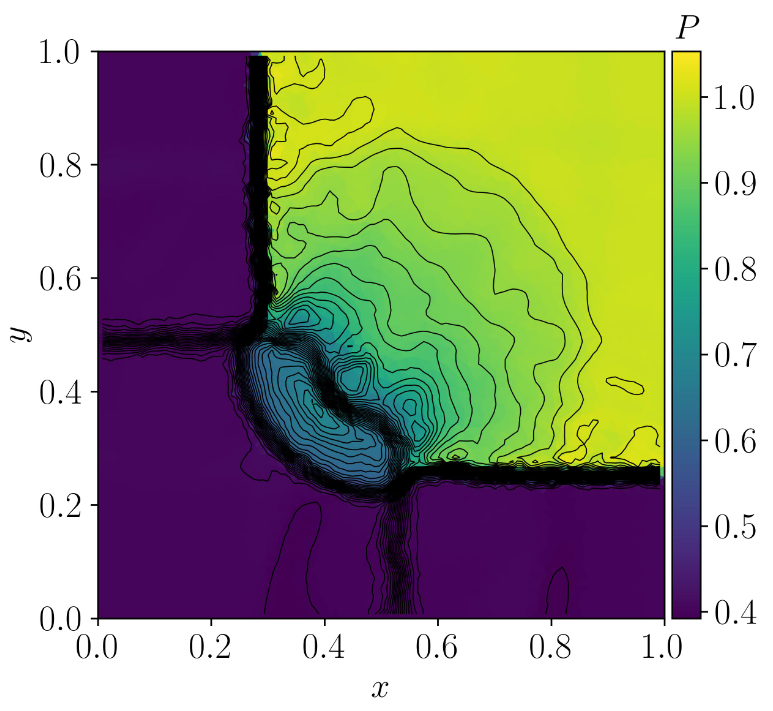}
         \caption{}
         \label{fig:sod_rho}
     \end{subfigure}
     \begin{subfigure}[b]{0.4\textwidth}
         \centering
         \includegraphics[width=\textwidth]{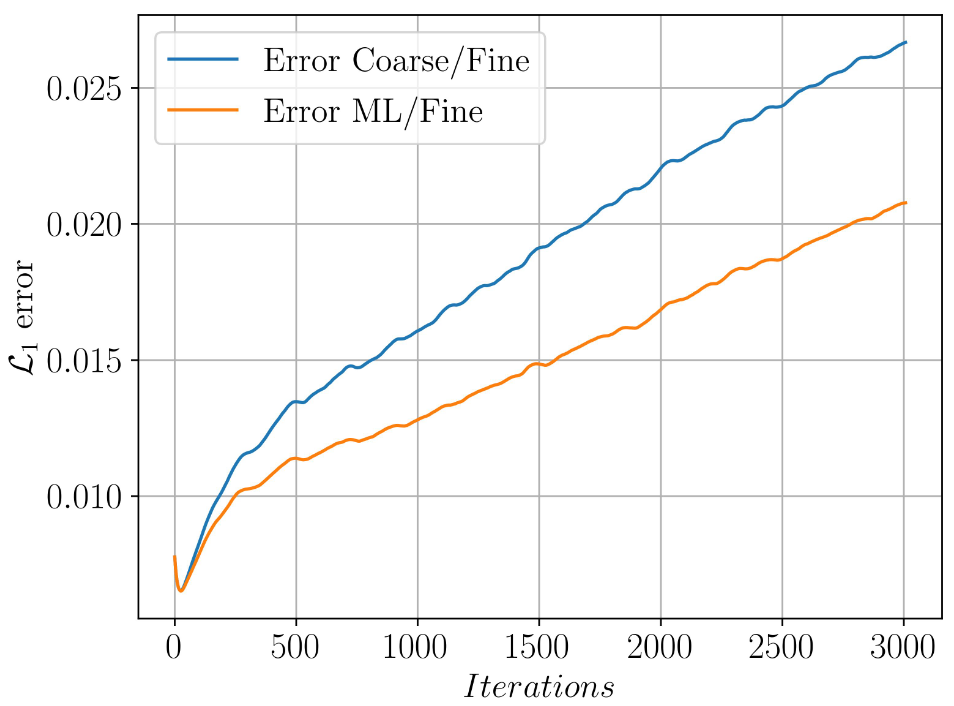}
         \caption{}
         \label{fig:sod_rho_zoom}
     \end{subfigure}
     \caption{Results for the 2D Riemann problem case 11 \citep{lax1998solution}. Density color map is overlayed by 30 density contours. The computations were performed on the  $(x, y) \in (0, 1) \times (0, 1)$ square with 5150 cells. (a) Pressure $p$ , (b) Error across iterations.}
     \label{fig:test_case_11}
\end{figure}

The test case 11 is made of contact waves and shocks. All boundary conditions have been set to subsonic outflow. It has a diagonal symmetry and the symmetry is preserved as seen on Figure \ref{fig:test_case_11}. The gain ($\sim 22\%$) is once again notably high and the front shocks are really straight. 

\begin{figure}
     \centering
     \includegraphics[width=0.45\textwidth]{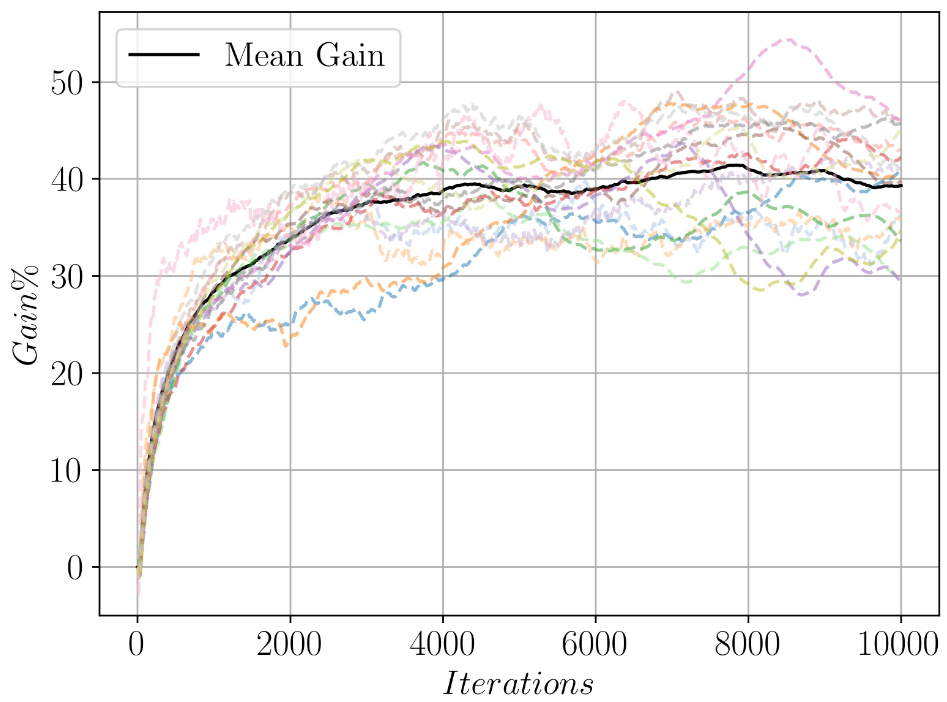}
     \caption{All gains across iterations for the 18 test cases present in \cite{lax1998solution} with the LSQ based gradient (Eq. (\ref{eq:ML_LSQ}). Each different dashed lines represent a different test case. The computations were performed on the  $(x, y) \in (0, 1) \times (0, 1)$ square with 5150 cells.}
     \label{fig:allCases}
\end{figure}

\begin{figure}
     \centering
     \includegraphics[width=0.45\textwidth]{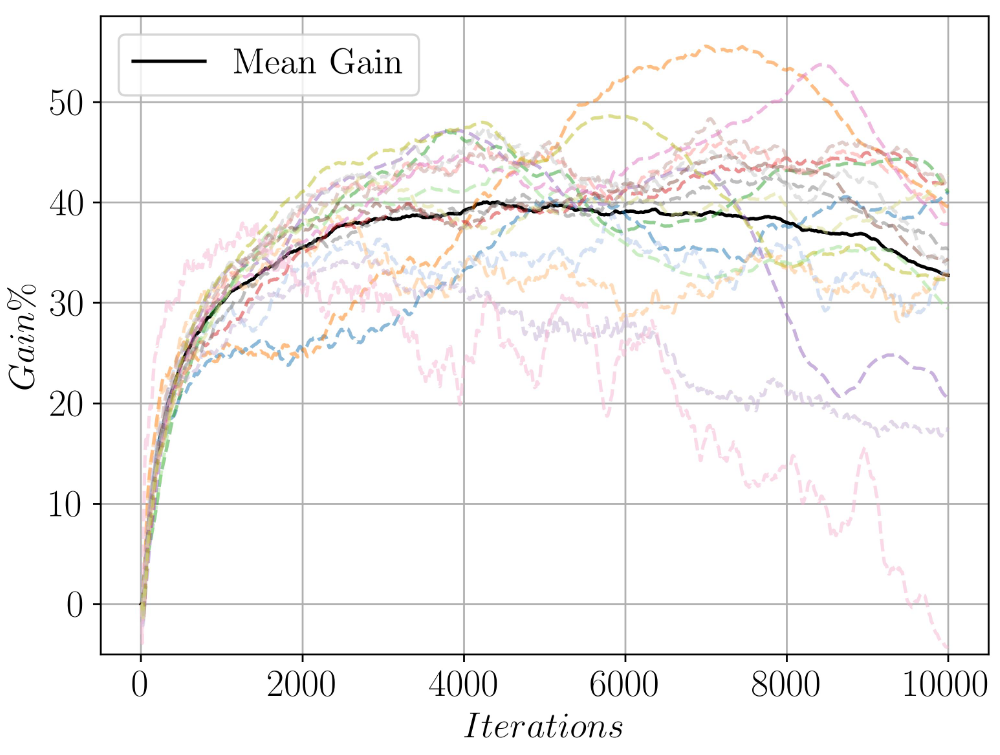}
     \caption{All gains across iterations for the 18 test cases present in \cite{lax1998solution} with the Green-Gauss based gradient (Eq. (\ref{eq:ML_GG}). Each different dashed lines represent a different test case. The computations were performed on the  $(x, y) \in (0, 1) \times (0, 1)$ square with 5150 cells.}
     \label{fig:allCases_GG}
\end{figure}

\paragraph{All test Riemann cases of \citep{lax1998solution}} In order to compare efficiently the method, all test cases from \cite{lax1998solution} have been computed to obtain a mean gain for 10000 iterations. The boundary conditions have been set to periodic for easier computation. As seen on Figure \ref{fig:allCases}, all test cases have a really good gain compared to the traditional finite volume solver. Just to note that some test cases have a diminishing gain after 10000 iterations step but the overall gain reaches 40\%.
To compare to the least squares based gradient computation, the same graph has been made on Figure \ref{fig:allCases_GG} for the Green-Gauss based gradient. Test cases with high values as initial conditions like test case 14 (pink dashed line reaching 0\% gain at 10000 time-steps on Figure \ref{fig:allCases_GG}) have really deceiving results. In general, results for test cases with symmetries lack symmetry for extended computation time. Moreover, the results are in general less precise probably because the Green-Gauss gradient is not first-order accurate unlike the LSQ gradient.

\subsection{The forward facing step}
The forward-facing step test case described by \cite{woodward1984numerical} is a challenging test case, mainly due to the rarefaction taking place at the corner of the step. This test case features multiple shock fronts across the domain and various boundary conditions (supersonic inflow/outflow, wall slip). This makes it an ideal test case for benchmarking. Therefore, we approach it as a two-dimensional problem. The initial condition reads $(\rho, u, v, p) = (1.4, 3., 0., 1.)$. The boundary conditions are a supersonic inflow on the left boundary of the domain, a supersonic outflow on the left boundary of the domain and the rest of the boundary conditions are slip-wall boundaries (top and bottom boundaries including the step).
\begin{figure}
     \centering
     \begin{subfigure}[b]{0.6\textwidth}
         \centering
         \includegraphics[width=\textwidth]{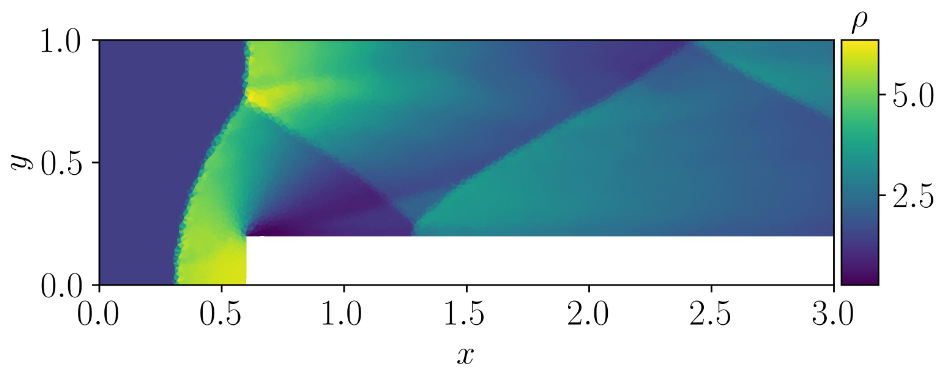}
         \caption{}
         \label{fig:forward_rho}
     \end{subfigure}
     \hfill
     \begin{subfigure}[b]{0.35\textwidth}
         \centering
         \includegraphics[width=\textwidth]{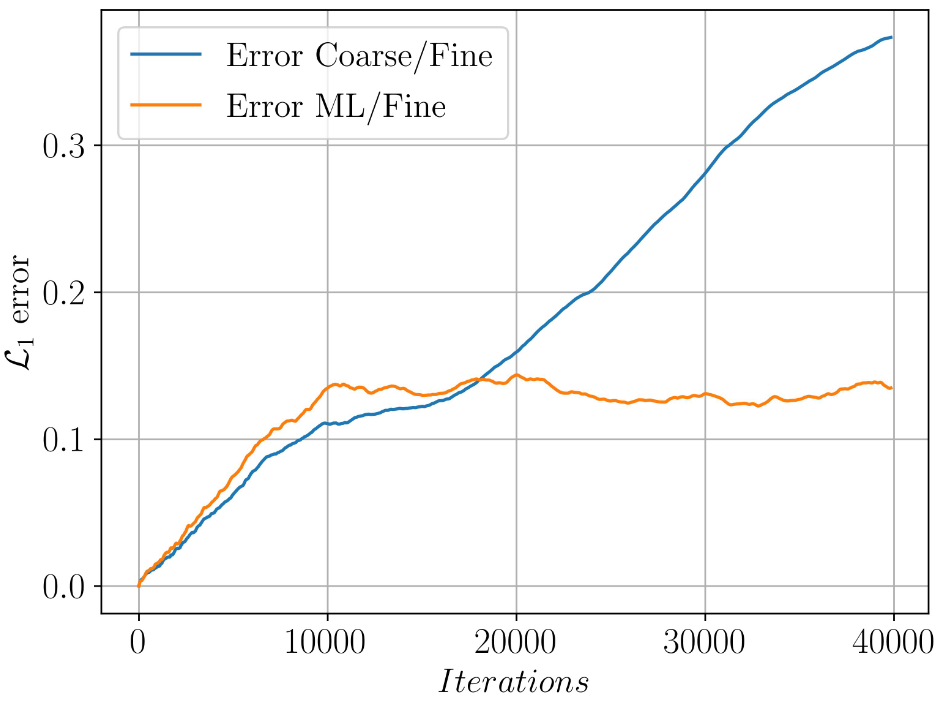}
         \caption{}
         \label{fig:forward_error}
     \end{subfigure}
     \caption{Results for the 2D forward facing step\cite{woodward1984numerical} at $T = 4$. the density solution field is represented. The computations were performed on a domain with 13089 cells. (a) Density $\rho$ , (b) Error across iterations. }
     \label{fig:forward}
\end{figure}
\par
\begin{figure}[H]
     \centering
     \includegraphics[width=0.45\textwidth]{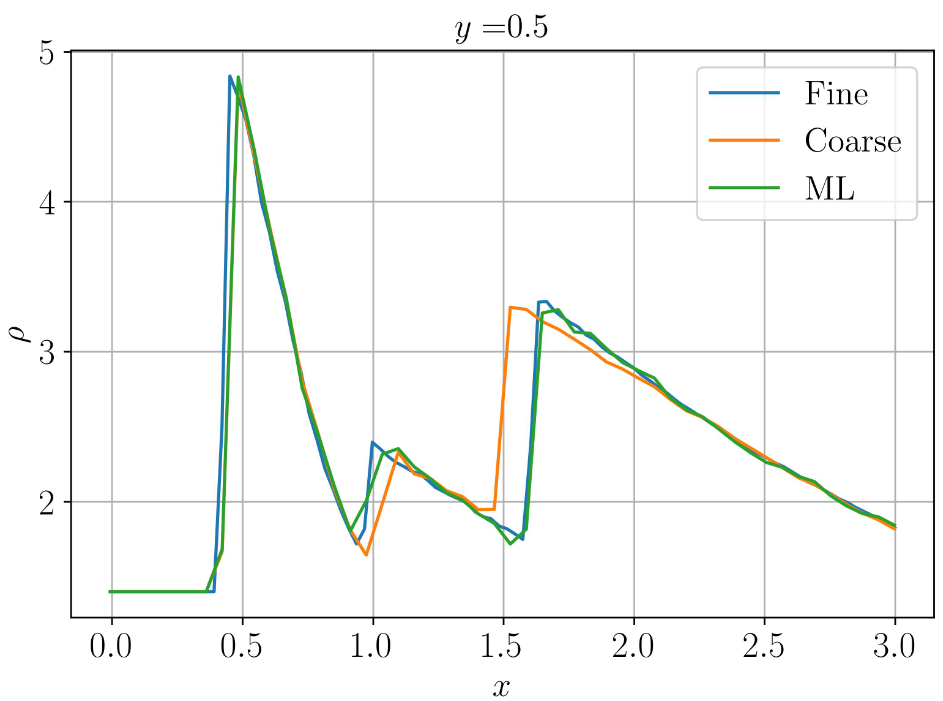}
     \caption{Slice at $y=0.5$ of density of the solution for the forward facing step represented Figure \ref{fig:forward}.}
     \label{fig:forward_slice}
\end{figure}

As seen on Figure \ref{fig:forward}, all shocks front are really well resolved, there seems to be no spurious oscillations. The resulting gain of the calculation is slightly over 60\%, although it was initially negative. For further investigation, we can note on Figure \ref{fig:forward_slice} that as opposed to the the traditional finite volume solver, the ML solver main advantage is its ability to accurately predict the shocks front on a coarse mesh.

\section{Convergence study} \label{sec:convergence}
\begin{figure}[H]
     \centering
     \includegraphics[width=0.45\textwidth]{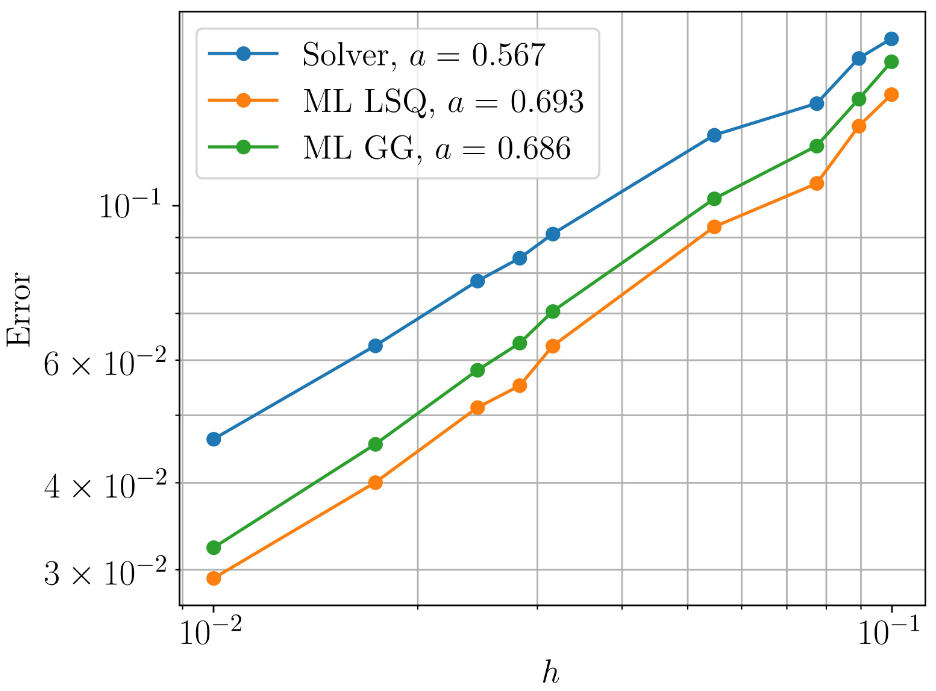}
     \caption{Mesh convergence with all Riemann initial conditions present in \cite{lax1998solution} computed for $T=0.2$ with periodic boundary conditions on a log-log scale. The traditional finite solver mesh convergence is represented in blue with a slope $a=0.567$, the ML with LSQ inspired gradient finite solver is represented in orange with a slope $a=0.693$ and the ML with GG inspired gradient finite solver is represented in green with a slope $a=0.686$. $h$ represents the average length of the cells present in the domain.}
     \label{fig:convergence}
\end{figure}

Figure \ref{fig:convergence} represents mesh convergence for the reference traditional second order finite volume solver and the mesh convergence for the ML finite volume solver. The convergence study has been conveyed using all Riemman test cases of \citep{lax1998solution} with periodic boundary conditions. 
The slope being inferior to the theoretical slope (2) obtained for the traditional solver is explained by the presence of strong shocks in the solutions used to compute the errors. 
The convergence slopes obtained for the GG based gradient and the LSQ based gradients are similar but as expected, the LSQ method is more precise. 

\section{Computational performance}  \label{sec:Computational}
\begin{figure}[H]
     \centering
     \includegraphics[width=0.45\textwidth]{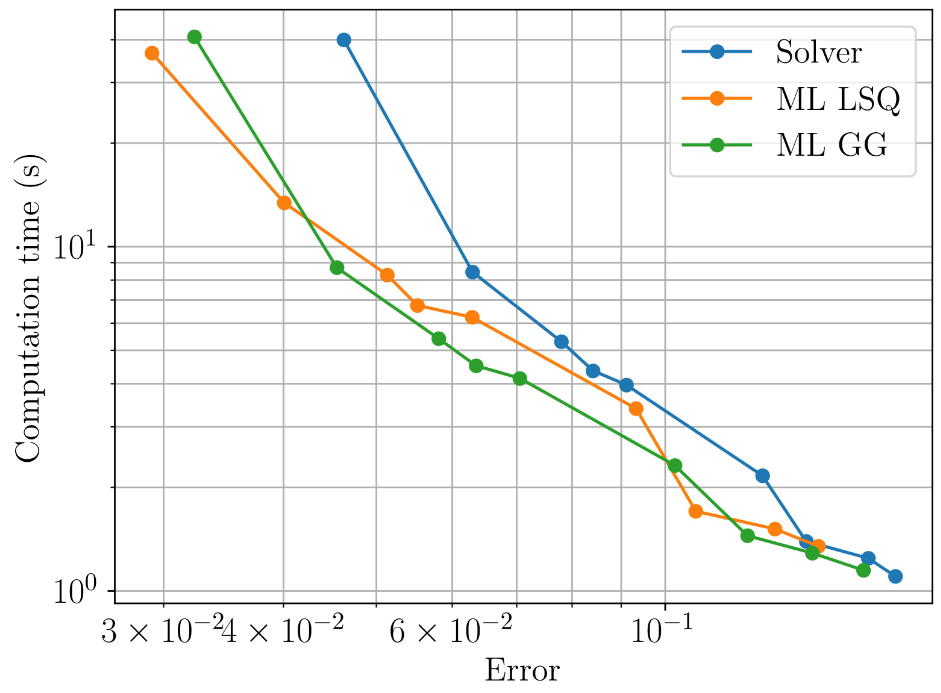}
     \caption{Computation time vs error on a log-log scale. The points used are the same as on Figure \ref{fig:convergence}.}
     \label{fig:error_vs_time}
\end{figure}

All points present on Figure \ref{fig:error_vs_time} has been obtained during the convergence study Sec. \ref{sec:convergence}. An error is then linked intrinsically to a cell size $h$. All calculations have been computed on exclusive GPU nodes for reproductibility.
\par
Figure \ref{fig:error_vs_time} represents the computation time for the ML solver and the traditional solver as the function of the error. As the error diminishes, the computation gain increases.
For coarse meshes (as illustrated on the right side of the figure), the time gain approaches zero. It is interesting to note that the computation time for a given mesh is shorter when using the ML solver as compared to the traditional solver. This suggests that the ML gradient requires less computation time as compared to the traditional LSQ gradient, which seems counterintuitive given the increased complexity of the ML gradient.
The complexity of the GG gradient is lower than the LSQ gradient and it allows better computational performances for a range of error. But as the error diminishes, the LSQ based ML gradient method becomes more efficient for a given error. 
\par 
We have to bear in mind that Figure \ref{fig:error_vs_time} does not take into account training time.

\section{Concluding remarks}
The work extended the previous work of \cite{de2024data} for unstructured finite volume solvers, focusing on the 2D Euler equations. The core innovation lies in extending the DeepONet architecture in order to build a local geometry aware gradient operator to have a better reconstruction of variables at the interfaces. To ensure the physical feasibility such as  entropy stability or total variation diminishing of the neural network solution to the forward problems, several regularizers have been introduced. The proposed method demonstrates notable improvements in solution accuracy (20–60\%) and computational efficiency, with successful generalization across benchmark problems and complex flow configurations. The least-squares based gradient formulation particularly showed robustness and higher convergence rates compared to traditional solvers.
\par
But the method doesn't come without some limitations. The trained model is somewhat dependent on specific mesh configurations, particularly the geometry used for training. Significant variations in mesh quality, anisotropy, or refinement strategies reduce accuracy or require retraining. A method exploring distances awareness can be developed for this case and can be applied for mesh-refinement. Moreover, performances degrades with extreme shocks strengths like for the double Mach reflection test case \cite{woodward1984numerical}. 
\\
The model has been trained using a specific CFL like condition $C_o$ for the time discretization scheme, one can change this condition but for large time steps the gain can become negative. Indeed, the ML method is slightly less robust regarding the CFL condition.
\\
During training, the quality of the resulting model is dependent on the the entropy and TVD penalizations used, they are heuristic and problem-specific. Furthermore, it requires an accurate fine-tuning of the parameters to build an efficient model.
\par 
Finally, for industrial purposes, one of the main limitation is the need of a fully differentiable solver for training even if some frameworks have been developed in low-level languages. One of the other drawbacks of the method for industrial applications is the need of a hybrid CPU-GPU architecture for efficient computing.

\section{Acknowledgments}
We thank Piotr Mirowski from Google DeepMind for helpful discussions regarding the model architecture. 

\bibliographystyle{elsarticle-harv} 
\bibliography{bibliography}

\newpage
\appendix
\section*{Appendices}
\addcontentsline{toc}{section}{Appendices}
\renewcommand{\thesubsection}{\Alph{subsection}}
\subsection{Dataset} \label{sec:dataset}
In this section are presented the initial conditions used for training on the domain $\Omega = [0,1]^2$.
\par
We first define the function 
\begin{equation*}
    \mathcal{C}(\textbf{x},\mathbf{p}) =
    \begin{cases}
        p_0 & \text{if $||(\textbf{x} - (0.5, 0.5)||_ 2\leq 0.125$}\\
        p_1 & \text{else if $x<0.5$ and $y<0.5$} \\
        p_2 & \text{else if $x\geq0.5$ and $y<0.5$} \\
        p_3 & \text{else if $x<0.5$ and $y\geq0.5$} \\
        p_4 & \text{else}
    \end{cases}
\end{equation*}
and the function
\begin{equation*}
    \mathcal{R}(\textbf{x},\mathbf{p}) =
    \begin{cases}
        p_0 & \text{else if $x<0.5$ and $y<0.5$} \\
        p_1 & \text{else if $x\geq0.5$ and $y<0.5$} \\
        p_2 & \text{else if $x<0.5$ and $y\geq0.5$} \\
        p_3 & \text{else if $x\geq0.5$ and $y\geq0.5$}
    \end{cases}
\end{equation*}
with $\textbf{p}$ being parameters following a distribution $\mathcal{U}(0,1)$. 
\par
Using these predefined functions, we define three types of functions used for initial conditions. Setting $\operatorname{max}$ as a parameter, we have
\begin{align*}
\textbf{f}_1(\textbf{x}) &=
    \begin{cases}
        f_{\rho} &= 0.5\operatorname{max} \times [a_0 sin(4 \pi x + \phi_0\pi) + a_1 sin(4 \pi y + \phi_1\pi)\\ 
        & + a_0 + a_1 + 0.1] \\ 
        f_{u} &= 3 [a_2 sin(4 \pi x + \phi_0\pi) + a_3 sin(4 \pi y + \phi_1\pi)]\\
        f_{v} &= 3 [a_4 sin(4 \pi x + \phi_0\pi) + a_5 sin(4 \pi y + \phi_1\pi)]\\
        f_{p} &= 0.5\operatorname{max} \times [a_6 sin(4 \pi x + \phi_0\pi) + a_7 sin(4 \pi y + \phi_1\pi)\\ 
        & + a_6 + a_7 + 0.1] 
    \end{cases} \\
\textbf{f}_2(\textbf{x}) &=
    \begin{cases}
        f_{\rho} &= \mathcal{C}(\textbf{x},\operatorname{max} \times \mathbf{p}_0 + \frac{1}{2}) \\ 
        f_{u} &= \mathcal{C}(\textbf{x},\operatorname{max} \times\mathbf{p}_1)\\
        f_{v} &= \mathcal{C}(\textbf{x},\operatorname{max} \times\mathbf{p}_2)\\
        f_{p} &= \mathcal{C}(\textbf{x},\operatorname{max} \times\mathbf{p}_3 + 0.2)
    \end{cases} \\
\textbf{f}_3(\textbf{x}) &=
    \begin{cases}
        f_{\rho} &= \mathcal{R}(\textbf{x},\operatorname{max} \times \mathbf{p}_0 + \frac{1}{2}) \\ 
        f_{u} &= \mathcal{R}(\textbf{x},\operatorname{max} \times\mathbf{p}_1)\\
        f_{v} &= \mathcal{R}(\textbf{x},\operatorname{max} \times\mathbf{p}_2)\\
        f_{p} &= \mathcal{R}(\textbf{x},\operatorname{max} \times\mathbf{p}_3 + 0.2)
    \end{cases}
\end{align*}
max is usually set to 6. All parameters $a_i$ follows a distribution $\mathcal{U}(0,1)$. The training data is composed of 8 initial conditions composed of 50\% $\textbf{f}_1$, 25\% $\textbf{f}_2$ and 25\% $\textbf{f}_3$  integrated on 2,000 time steps.
\begin{figure}[ht]
\captionsetup{justification=centering}
    \begin{subfigure}[h]{0.3\textwidth}
        \includegraphics[width=\textwidth]{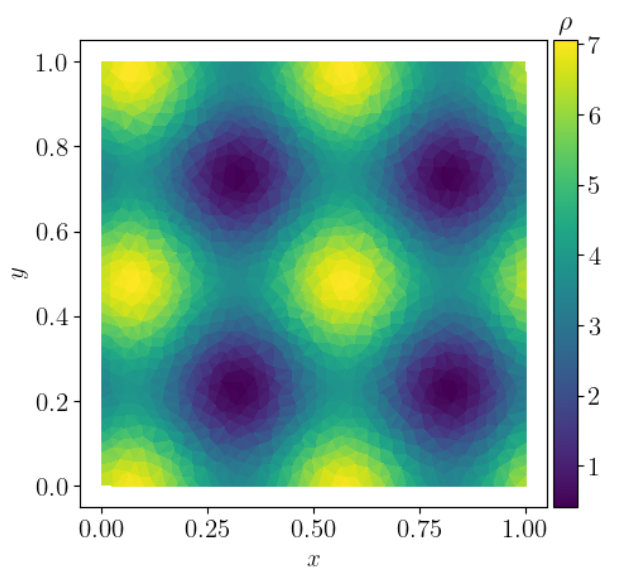}
        \caption{A representation of $\textbf{f}_1$}
    \end{subfigure}
    \begin{subfigure}[h]{0.3\textwidth}
        \includegraphics[width=\textwidth]{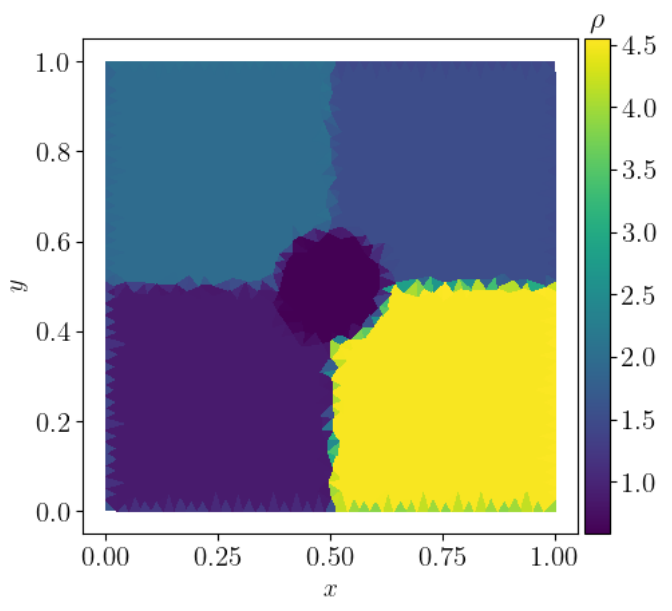}
        \caption{A representation of $\textbf{f}_2$}
    \end{subfigure}
    \begin{subfigure}[h]{0.3\textwidth}
        \includegraphics[width=\textwidth]{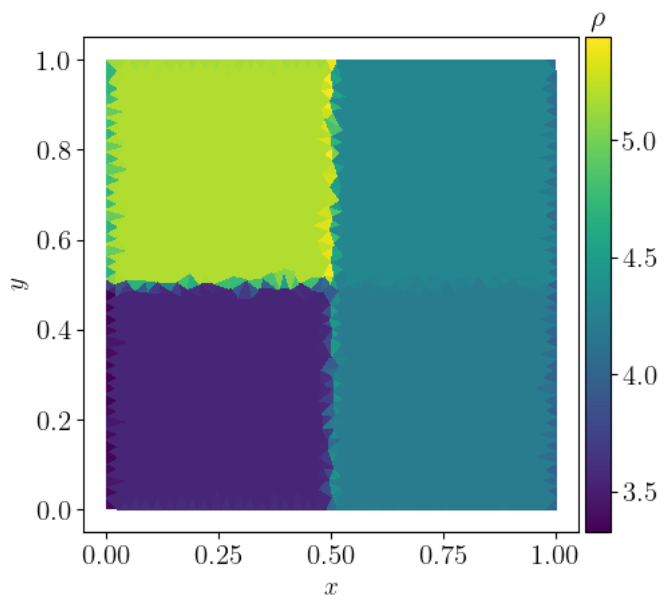}
        \caption{A representation of $\textbf{f}_3$}
    \end{subfigure}
\caption{Initial conditions for the 2D Euler database with periodic boundary conditions.}
\label{fig:database_2D_euler}
\end{figure}

\end{document}